\documentclass[11pt]{article}
\usepackage{amsmath}
\usepackage{amsfonts}
\usepackage{amsthm}
\usepackage{amssymb, setspace}
\usepackage{color,graphicx,epsfig,geometry,hyperref,fancyhdr}
\usepackage[T1]{fontenc}
\usepackage{float}
\usepackage{subfig}
\usepackage{commath}
\usepackage{bbm}
\usepackage{mathrsfs,fleqn}
\usepackage{pgfplots}
\pgfplotsset{compat=newest}

\usepackage{multirow}

\usepackage[pagewise]{lineno} 

\begin{document}

\title{Application of machine learning regression models to inverse eigenvalue problems}

\date{}
\maketitle

\centerline{\scshape Nikolaos Pallikarakis \footnote{Corresponding author}}
\medskip
{\footnotesize
 \centerline{Department of Mathematics, National Technical University of Athens}
  \centerline{ 157 80, Athens, GREECE}
  
 \centerline{n\&k Technology, Inc.}
  \centerline{CA 95119, San Jose, USA}
  \centerline{Email:  \texttt{npall@central.ntua.gr}}

} 

\medskip

\centerline{\scshape Andreas Ntargaras}
\medskip
{\footnotesize
 \centerline{n\&k Technology, Inc.}
    \centerline{CA 95119, San Jose, USA}
    \centerline{Email:  \texttt{andreas@nandk.com}}
}


\begin{abstract}
In this work, we study the numerical solution of inverse eigenvalue problems from a machine learning perspective. Two different problems are considered: the inverse Strum-Liouville eigenvalue problem for symmetric potentials and the inverse transmission eigenvalue problem for spherically symmetric refractive indices. Firstly, we solve the corresponding direct problems to produce the required eigenvalues datasets in order to train the machine learning algorithms. Next, we consider several examples of inverse problems and compare the performance of each model to predict the unknown potentials and refractive indices respectively, from a given small set of the lowest eigenvalues. The supervised regression models we use are k-Nearest Neighbours, Random Forests and Multi-Layer Perceptron. Our experiments show that these machine learning methods, under appropriate tuning on their parameters, can numerically solve the examined inverse eigenvalue problems.\\ 
\textbf{\textit{keywords:}} inverse spectral problems, transmission eigenvalues, machine learning, deep learning

\end{abstract}

\section{Introduction}

Application of machine learning and deep learning algorithms to solve inverse problems is a very active topic the latest years. The research in this area has been focused on imaging and tomography (see e.g.\cite{JMRU, LSAH, OJMBD} and references therein), on nonlinear differential equations of mathematical physics \cite{RPK}, on quantum mechanics and physics \cite{CXZHZZ,CAFL} and other topics of applied sciences. Among other models, deep neural networks are one of the state-of-the art methods focused on the ill-posed nature of the problems \cite{AO,AMOS, LSAH}. 

An interesting question is whether machine learning can be also applied to solve inverse eigenvalue problems. That is, to train an algorithm using spectral data and examine it's ability to predict an unknown function, which is the coefficient of a differential operator and is associated with the physical properties of the problem under consideration. Solution of direct eigenvalue problems has been considered several years earlier \cite{CU} and is ongoing. See for example \cite{KEKFHGBOYS} for an application of neural networks to solve eigenvalue problems in micromagnetics. However, the idea of using learning models to deal with the inverse problem is more recent. A result on the inverse problem for the Schr\"odinger operator in one dimension using neural networks, is given in \cite{LACP}. Furthermore, we refer to \cite{ORR} for an application of neural networks for the inverse eigenvalue problem of the anisotropic Dirichlet operator and to \cite{OR} for the estimation of the Lamé coefficients from the eigenvalues of the elasticity operator. In these last two papers, authors applied finite element methods to obtain the eigenvalues, which where afterwards used to train the corresponding algorithms. 

In the present work, we apply supervised machine learning regression models, which means that the algorithms are trained given a set of input or feature values (the eigenvalues) together with a set of output or target values (the coefficients of the differential operator). 
Our goal is to show that given a sufficient sample of spectral data, machine learning algorithms are capable to get trained and predict the unknown functions for the solution of inverse problems. This approach, could be potentially usefully to the numerical solution of problems of this type together (and not necessarily in contrast) with more classical approaches, like optimization.  

We study two different eigenvalue problems:  the self-adjoint Strum-Liouville eigenvalue problem for symmetric potentials, as a toy problem, and the non-self-adjoint transmission eigenvalue problem with spherically symmetric piecewise constant refractive index. Firstly, we formulate these eigenvalue problems and pose the numerical solution of the corresponding direct problems. Next, we present the basic properties of the machine learning algorithms we implement, which are k-Nearest Neighbours (kNN), Random Forests (RF) and Multi-Layer Perceptron (MLP), as well as the data pre- and post-processing methods we use. Afterwards, we numerically solve the direct problems for various examples and create the spectral data which in turn are used as training data for the machine learning models. The direct Strum-Liouville eigenvalue problem is solved using the matlab package \textsc{matslice} \cite{LVD} and the direct transmission eigenvalue problem by utilizing a spectral-Galerkin method introduced in \cite{GP}. Finally, we present numerical results for the inverse problems and compare the efficiency of each algorithm, to train and predict. Testing of each model is carried out with out-of-sample data, i.e. eigenvalues that are not included in the training/validating sets, but fall within the range of these sets. For the inverse transmission problem, this is achieved using a few of the lowest eigenvalues calculated from separation of variables. Each model is evaluated using the coefficient of determination $R^2$ and the root mean square error, where we also evaluate the importance of each feature (eigenvalue) to train the models, using an input perturbation algorithm. To our knowledge, this is a first attempt on the numerical solution of the inverse transmission eigenvalue problem, from a machine learning perspective. 

\section{Formulation of the spectral problems} \label{form}

\subsection{Sturm–Liouville eigenvalue problem} 
We consider the Sturm-Liouville eigenvalue problem in canonical form: 
\begin{equation}
    -y^{\prime\prime}+q(x)y=\lambda y, \qquad 0<x<1, \label{sl}
\end{equation}
with boundary conditions: 
\begin{equation}
   y^{\prime}(0)-hy(0)=0, \qquad   y^{\prime}(1)+Hy(1)=0, \label{slbc}
\end{equation}
where the potential function $q(x)$ is in $L^1(0,1)$. The values of $\lambda$ for non-trivial solutions $y$, are the eigenvalues of the problem. The above is one of the most studied self-adjoint boundary value problem for which there exists a comprehensive spectral theory concerning both the direct and the inverse problem  \cite{FY, Mar}. 

The direct spectral problem consists of finding the eigenpairs $(\lambda;y)$ for a given potential $q$. On the other hand, the inverse problem is to recover the unknown potential using given spectral data. For the special case of a symmetric potential $q(x)=q(1-x)$ with $h=H$, it is a well known and old result that one spectrum $\{\lambda_l\}_{l=1}^{\infty}$ suffices to uniquely determine $q$ \cite{Bo}.

\subsubsection{Numerical solution of the direct problem} 
There have been several numerical methods developed to solve the direct problem (\ref{sl})-(\ref{slbc}) (see e.g. \cite{AMR,Py}). For the purpose of this work, we restrict ourselves to the case of a symmetric potential. We use the matlab package \textsc{matslice}, which implements piecewise constant perturbation methods to solve the boundary value problems, and compute the eigenvalues for a family of given potentials. With this software, we are able to calculate as many eigenvalues as needed and create the required data for training and validating the machine learning algorithms. 

\subsection{Transmission eigenvalue problem} 

We examine the interior transmission problem corresponding to the acoustic scattering problem of an isotropic inhomogeneous medium $D\subseteq \mathbb{R}^{2,3}$, with a smooth boundary $\partial D$ \cite{CK}. 

The interior transmission eigenvalue problem is a non-self-adjoint problem, defined as follows: Find $k\!\in\!\mathbb{C}$ and a non-trivial solution $w,v\in\!L^2(D)$ such that $w-v\!\in\!H_0^2(D)$ satisfying:
\begin{eqnarray}
&\Delta w+k^2n(x)w=0&\ \textrm{in}\ D, \label{trB1}\\
&\Delta v+k^2v=0&\ \textrm{in}\ D, \label{trB2}\\
&w=v&\ \textrm{on}\ \partial D, \label{trB3}\\
&\frac{\partial w}{\partial \nu}=\frac{\partial v}{\partial \nu}&\ \textrm{on}\ \partial D. \label{trB4}
\end{eqnarray} 

The complex numbers $k$ are the transmission eigenvalues and pairs $(w, v)$ are the corresponding eigenfunctions of the problem. The refractive index $n(x)$ is assumed to be a real, positive function on $L^{\infty}(D)$, not identically equal to $1$. It is shown in \cite{CGH} that there exists an infinite set of real transmission eigenvalues accumulating at $+\infty$, provided that either $0<n(x)<1$ or $n(x)>1$.

Since real transmission eigenvalues can been measured from scattering data \cite{CCH2} and are also associated with the material properties of the non-homogeneous medium \cite{CCaC2,CCM,Ha},  the transmission eigenvalue problem is a very active topic in inverse scattering theory (we refer to \cite{CCH,CH} and the references therein). Research on transmission eigenvalues has been focused both on theoretical and numerical aspects, in the direction of the direct and the inverse problem as well. 

In the direct spectral problem, we assume that the refractive index is known and we want to estimate the eigenpairs $\{k;(w,v)\}$. For the inverse spectral problem, given information on the spectrum and/or the eigenfunctions, we intend to recover the properties of the unknown refractive index. 

In the case when the domain is a ball of $\mathbb{R}^{2,3}$ and the refractive index $n(|x|):=n(r)$ is a function depending only on the radius, the method of separation of variables can be applied and (\ref{trB1})-(\ref{trB4}) is simplified to a non-self-adjoint boundary-value problem for an ordinary differential equation, with the spectral parameter $k$ appearing in the boundary condition, at the right endpoint \cite{MP}. Uniqueness for the refractive index $n(x)$ from the knowledge of all eigenvalues is shown in \cite{CCG} for a $C^2$ refractive index and in \cite{GPd} for a piecewise $C^2$ refractive index, under the assumption of either $n(x)>1$ or $0<n(x)<1$. We refer to \cite{AGP,BYY,CPS} for some further results with regards to existence and discreteness of eigenvalues and uniqueness of the inverse problem, when the domain is spherically symmetric.

\subsubsection{Numerical solution of the direct problem for a piecewise constant refractive index } \label{teppw}
We consider the simple case where the refractive index is a piecewise constant function and $D$ is the unit disc of $\mathbb{R}^2$. We assume that the domain has $L$-layers such that $D=\cup_{i=1}^{L}D_i$ and $\partial D_i$ are concentric circles. An example is shown in the following figure,  
\begin{figure}[hbtp]
\centering
\includegraphics[scale=0.55]{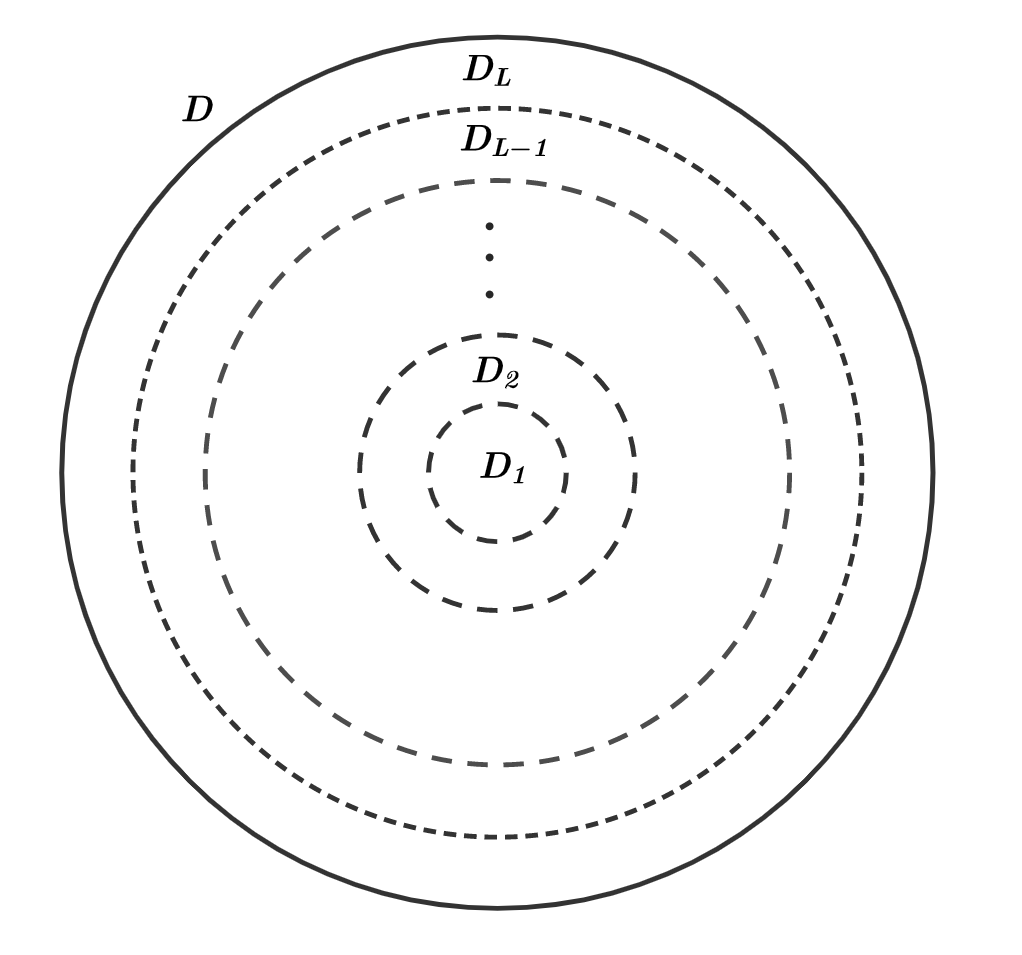}
\caption{A disc with $L$-layers}
\end{figure}
where the piecewise constant index of refraction $n(r)$ has $L-1$ number of discontinuities at $\{d_l\}_{l=1}^{L-1}$ and is given by:
\begin{equation}
n(r)=\left\{
\begin{array}{ll}
n_1,& x\in D_1\\
\vdots\\
 n_L,& x\in D_L. \label{refr}
\end{array} \right.
\end{equation}
Refractive indices of this type can model for example waveguides with circular cross-section, like multi-step index fibres \cite{SL}.  

By applying the method of separation of variables, we derive that the eigenfunctions of (\ref{trB1})-(\ref{trB2}) have the following form:
\begin{equation}
\begin{aligned}
    v_m(r,\theta) = &\ a_{m}J_m(k r)e^{im\theta},\hspace{5.3cm} 0<r<1, \\
     w_m(r,\theta)= &\begin{cases} 
      b_{1,m}J_m(k\sqrt{n_1}r)e^{im\theta},& 0<r<d_1, \\
 \left(b_{2,m}J_m(k\sqrt{n_2}r)+c_{2,m}N_m(k\sqrt{n_2}r)\right)e^{im\theta}, & d_1<r<d_2,\\
      \vdots\\      \left(b_{L,m}J_m(k\sqrt{n_L}r)+c_{L,m}N_m(k\sqrt{n_L}r)\right)e^{im\theta}, & d_{L-1}<r<1,
      \end{cases}   
\end{aligned} \label{eignf}
\end{equation}
for $m \in \mathbb{Z}$, where $a_m,b_{m},c_{m}$ are 
 constants and  $J_m$ and $N_m$ are Bessel and Neumann functions, respectively. By solving the system of the boundary transmission conditions (\ref{trB3})-(\ref{trB4}) and the continuity of the eigenfunctions (\ref{eignf}) at the jumps $\{d_l\}_{l=1}^{L-1}$, we conclude that $k$ is a transmission eigenvalue, if and only if is a root of the determinant $D_m(k)$ of the following $2L\times 2L$ matrix: 
 
 \footnotesize{
\begin{equation}\label{detN}
\begin{aligned}
&\vspace{-2cm}\left[\begin{matrix}
 J_m\left(kr\right)|_{r=1}&0&0&0& \cdots\\
  \frac{\mathrm{d}}{\mathrm{d}r}J_m(k r)\big{|}_{r=1}&0&0&0& \cdots\\
  0&J_m\left(k\sqrt{n_1}r\right)|_{r=d_1}&-J_m\left(k\sqrt{n_2}r\right)|_{r=d_1}&-N_m\left(k\sqrt{n_2}r\right)|_{r=d_1}&\cdots\\
  0&\frac{\mathrm{d}}{\mathrm{d}r}J_m\left(k\sqrt{n_1}r\right)\big{|}_{r=d_1}&-\frac{\mathrm{d}}{\mathrm{d}r}J_m\left(k\sqrt{n_2}r\right)\big{|}_{r=d_1}&-\frac{\mathrm{d}}{\mathrm{d}r}N_m\left(k\sqrt{n_2}r\right)\big{|}_{r=d_1}&\cdots\\
 \vdots&\vdots&\vdots&\vdots&\vdots\\
    0&0&0&0&\cdots\\
    0&0&0&0&\cdots\\
\end{matrix}\right.\\
&\hspace{-1.2cm}
\left.\begin{matrix}
0&0&-J_m\left(k\sqrt{n_L}r\right)|_{r=1}&-N_m\left(k\sqrt{n_L}r\right)|_{r=1}\\
0&0&-\frac{\mathrm{d}}{\mathrm{d}r}J_m\left(k\sqrt{n_L}r\right)\big{|}_{r=1}&-\frac{\mathrm{d}}{\mathrm{d}r}N_m\left(k\sqrt{n_L}r\right)\big{|}_{r=1}\\
0&0&0&0\\
0&0&0&0\\
\vdots&\vdots&\vdots&\vdots\\
J_m\left(k\sqrt{n_{L-1}}r\right)|_{r=d_{L-1}}&N_m\left(k\sqrt{n_{L-1}}r\right)|_{r=d_{L-1}}&-J_m\left(k\sqrt{n_{L}}r\right)|_{r=d_{L-1}}&-N_m\left(k\sqrt{n_{L}}r\right)|_{r=d_{L-1}}\\
\frac{\mathrm{d}}{\mathrm{d}r}J_m\left(k\sqrt{n_{L-1}}r\right)\big{|}_{r=d_{L-1}}&\frac{\mathrm{d}}{\mathrm{d}r}N_m\left(k\sqrt{n_{L-1}}r\right)\big{|}_{r=d_{L-1}}&-\frac{\mathrm{d}}{\mathrm{d}r}J_m\left(k\sqrt{n_{L}}r\right)\big{|}_{r=d_{L-1}}&-\frac{\mathrm{d}}{\mathrm{d}r}N_m\left(k\sqrt{n_{L}}r\right)\big{|}_{r=d_{L-1}}
\end{matrix}\right]
\end{aligned}
\end{equation}}

\normalsize

\begin{figure}[ht] 
       \centering
    \includegraphics[scale=0.7]{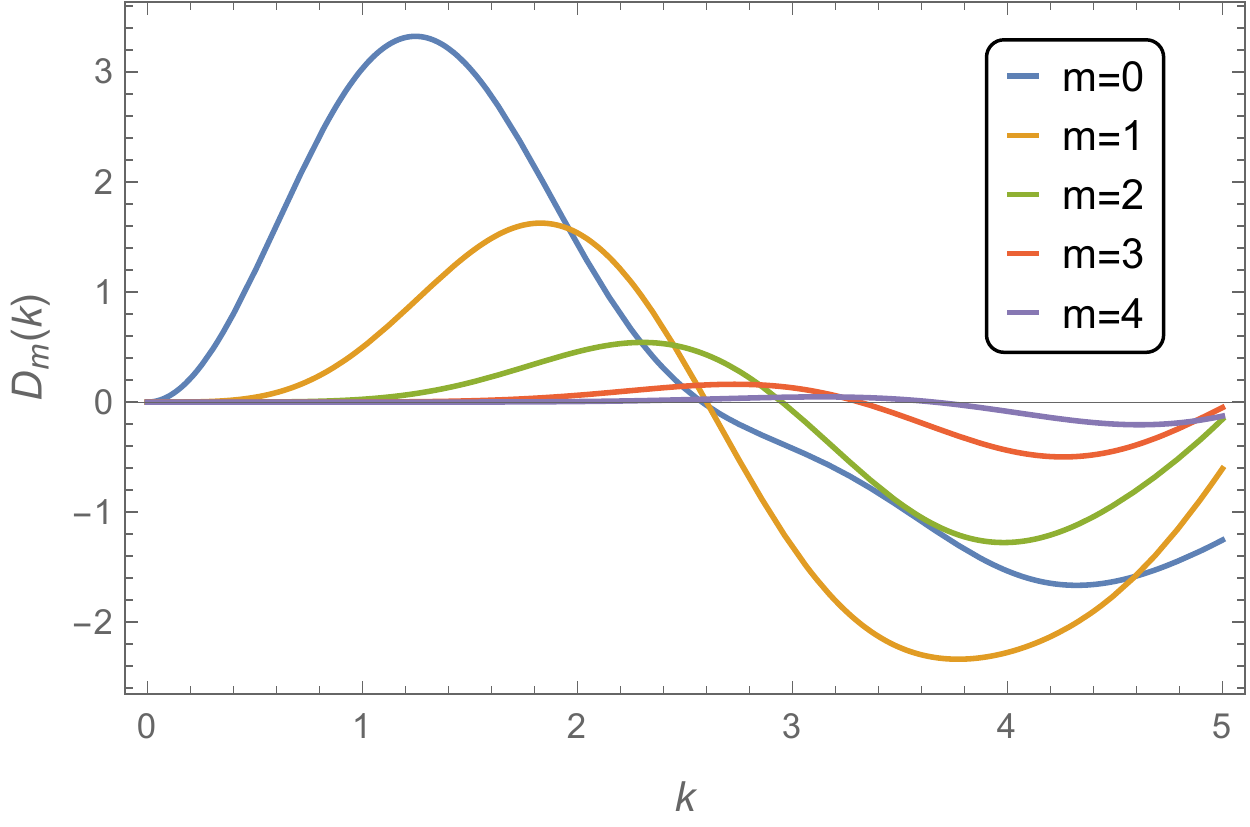} 
  \caption{Real transmission eigenvalues of the piecewise constant refractive index $(2.4, 7.4, 8.1, 2.1)$ with $d_1=0.25,\ d_2=0.5$ and $d_3=0.75$, for $m=0,1,2,3,4$.}
  \label{4L_det}
\end{figure}

We use this relation to compute the original eigenvalues with root-finding software, for discs with piecewise constant indices.  An example is shown in the figure \ref{4L_det} for a refractive index with $4$ layers.  

Several numerical methods to approximate the transmission eigenvalues in general domains have been proposed (we refer for example to \cite{An,CaMS,CoMS,Kl}). Here, we use the spectral-Galerkin method presented in \cite{GP}. The idea is to transform the problem (\ref{trB1})-(\ref{trB4}) to an equivalent fourth order eigenvalue problem for $u=w-v$ in the Sobolev space $H_0^2(D)$ and use the corresponding variational formulation: 
\begin{equation*}
\int_D\frac{1}{n-1}(\Delta u+k^2  u)(\Delta
\overline{\phi}+k^2n\overline{\phi})\mathrm{d}x=0,\ \ \forall\ \phi\in
H_0^2(D). 
\end{equation*}
By introducing a Hilbert basis $\{\phi_i\}_{i=1}^{N}$ of dimension $N$, the discrete Galerkin scheme for the above formula, is written in matrix form as:
\begin{equation}
[A^{(N)}-(k^{(N)})^2 B^{(N)}+(k^{(N)})^4C^{(N)}]\mathbf{c}=0, \label{a8}
\end{equation}
where 
\begin{eqnarray*}
&A^{(N)}:=&\int_{D}\frac{1}{n(x)-1}\Delta \phi _i\Delta \overline{\phi_j}\mathrm{d}x \\
&B^{(N)}:=&- \int_{D}\frac{n(x)}{n(x)-1}\Delta \phi _i\overline{\phi_j}\mathrm{d}x-\int_{D}\frac{1}{n(x)-1}\phi _i\Delta \overline{\phi_j}\mathrm{d}x \\
&C^{(N)}:=&\int_{D}\frac{n(x)}{n(x)-1} \phi _i \overline{\phi_j}\mathrm{d}x
\end{eqnarray*}
are $N\times N$ matrices, $\mathbf{c}=(c_1,c_2,\dots,c_N)^{\top},\ i,j=1,\dots,N$. Eigenvalues $k^{(N)}$ of the quadratic eigenvalue problem (\ref{a8}) are shown to converge to the corresponding eigenvalues $k$ of (\ref{trB1})-(\ref{trB4}), for $N\rightarrow\infty$ \cite[Prop. 2]{GP}. 


The above approximation method can be applied in domains with general geometry. Here, we restrict to spherically symmetric mediums with a refractive index of form (\ref{refr}). We can solve the direct transmission problem, using $2d$-numerical integration on $D$ and compute the corresponding matrices $A^{(N)}, B^{(N)} $ and $C^{(N)}$, where: 
\begin{eqnarray*}
&A^{(N)}:=&\sum_{l=1}^{L}\frac{1}{n_l-1}\int_{D_l}\Delta \phi _i\Delta \overline{\phi_j}\mathrm{d}x \\
&B^{(N)}:=&-\sum_{l=1}^{L}\frac{n_l}{n_l-1} \int_{D_l}\Delta \phi _i\overline{\phi_j}\mathrm{d}x-\sum_{l=1}^{L}\frac{1}{n_l-1}\int_{D_l}\phi _i\Delta \overline{\phi_j}\mathrm{d}x \\
&C^{(N)}:=&\sum_{l=1}^{L}\frac{n_l}{n_l-1}\int_{D_l} \phi _i \overline{\phi_j}\mathrm{d}x
\end{eqnarray*}
Then, the eigenvalues are calculated using the matlab function \textsc{polyeig}. With this method, we can directly compute the real
and complex transmission eigenvalues at the same time, with low computational cost depending on the size $N$ of the matrices. Within this approach, we generate the required spectral dataset which is used to train and validate the machine learning algorithms.

\section{Machine learning regression methods}

For the convenience of the reader, we present the basic properties of the machine learning models we implement in this work. Furthermore, we describe the pre-processing method we use as well as the regression metrics to estimate the efficiency of each algorithm. Finally, we present an algorithm which is used to evaluate feature importance. We refer to the textbooks \cite{GBC, HTF}, among others, for a more detailed reading.

\subsection{k-Nearest Neighbours} 

k-Nearest Neighbours (kNN) is a relatively simple supervised algorithm, firstly introduced in \cite{FH}. In the case of a regression model, the input consists of data $\mathrm{X}$ with the assumption of similarity among the k-closest of them and the output of values $\mathrm{y}$ which are the average of the k-closest neighbours. kNN learns the training data and predicts the value(s) for new outputs based on the k-nearest samples. Two important parameters of this model is the number of neighbours k and the distance metric with which the similarity is measured. The most commonly used metric is the euclidean distance, which for $N$-dimensional samples X and Z is: 
\begin{equation*}
    \textrm{dist} = \sqrt{(\mathrm{X}_1-\mathrm{Z}_1)^2+\cdots (\mathrm{Z}_N-\mathrm{Z}_N)^2}.
\end{equation*}
Given an input of new data $\mathrm{X}_{new}$, the predicted value of $\mathrm{y}_{pred}$ from k-closest neighbours $\mathrm{y}_i,\ i\in \mathrm{k}$ (estimated by $\mathrm{X}_{new}$), is given by the average: 
\begin{equation*}
    \mathrm{y}_{pred}(\mathrm{X}_{new}) = \frac{1}{\mathrm{k}}\sum_{i\in \mathrm{k}(\mathrm{X}_{new})}\mathrm{y}_i.
\end{equation*}
The number of neighbours used is very important. Small values of k may lead to overfitting (the algorithm is strongly attached to the training data and cannot generalize/predict in out-of-sample data) due to the result that is trained by using few samples. On the other hand, high values of k may include irrelevant data during the training process and as a result reduce the prediction accuracy. 

\subsection{Random Forests}

Random Forests (RF), is a supervised learning method firstly proposed  in \cite{Br}. It belongs in the category of ensemble methods, i.e. methods that use multiple different models trained in the same data and then average their outputs in an effort to to give a more accurate result. In RF, the different models are the decision trees, which are node -like models based on logical (if-else) conditions for splitting the nodes and ultimately creating the "forest".

\begin{figure}[ht] 
    \centering
    \includegraphics[scale=0.45]{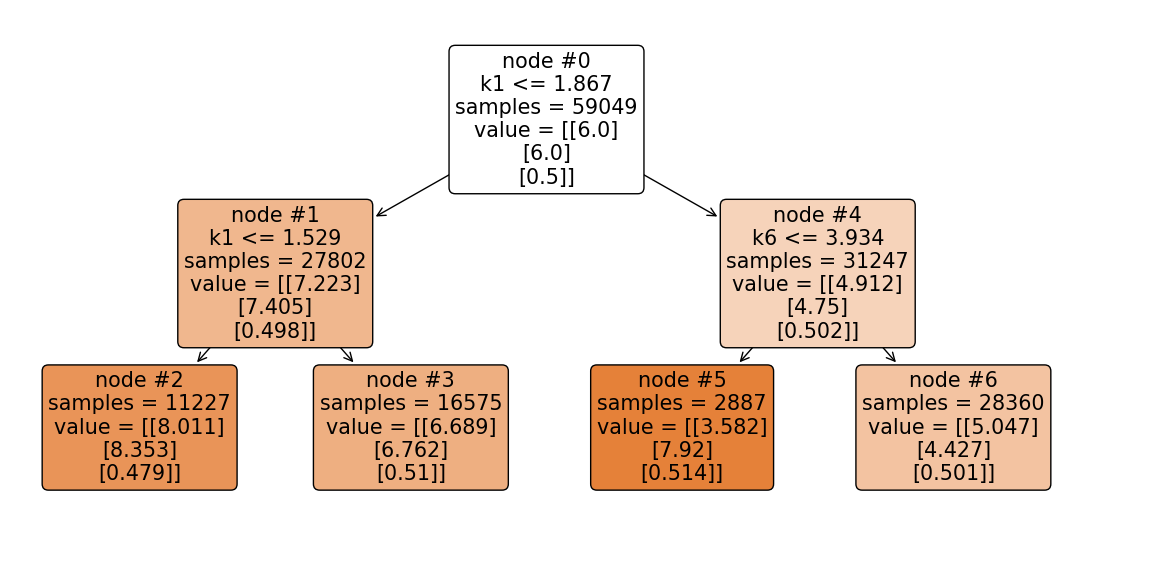} 
  \caption{A decision tree with depth $2$, corresponding to the training data of subsection \ref{sec421}.}
  \label{rf_tree}
\end{figure}

An example of a tree with depth equal to $2$ is shown in figure \ref{rf_tree}. This corresponds to the training data of the RF model used for the problem in subsection \ref{sec421}. The input values (transmission eigenvalues) are following an "if-else" inspection, starting from node$\# 0$ and ending in node$\#5$, notated with darker color, which corresponds to the predicted value (the refractive index).  

The steps for constructing a RF for a regression problem are the following: $(i)$ randomly choose $n$ samples of the training data set, $(ii)$ create a decision tree using these samples, $(iii)$ for a given value of $p$ trees, repeat steps $(i)$ and $(ii)$, $(iv)$ for any new input values, each tree does a prediction and the final output value will be the average of the predictions of all trees. 

 The most common criterion to determine the model outcomes is the mean square error. Other important parameters are the number of trees, the maximum depth of each tree as well as the number of features to consider when deciding for the best split. RF with a few trees are prone to overfit, so choosing the appropriate number of trees as well as their depth is essential. 

\subsection{Multi-Layer Perceptron} 

\begin{figure}[ht] 
    \centering
    \includegraphics[scale=0.35]{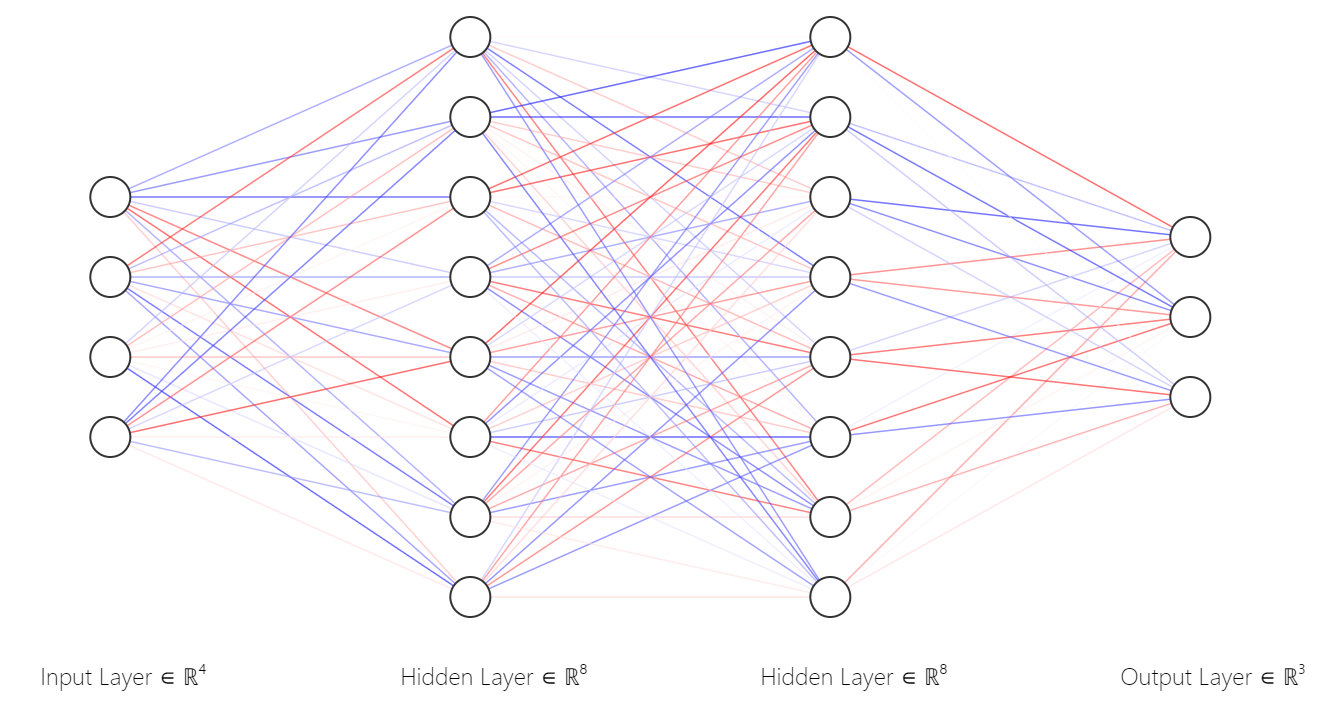} 
  \caption{A feed-forward ANN with two hidden layers \cite{LN}.}
  \label{ANN}
\end{figure}

Multi-Layer Perceptron (MLP) is an artificial neural network (ANN) architecture, belonging in the category of feed-forward networks. These are networks for which the connections between the layers doesn't form loops and the information is going only forward, from the input layer thought the hidden layer(s), to the output layer. An example is shown in figure \ref{ANN}. Each layer consists of a number of nodes (neurons) and each node is connected with a specific weight, with the nodes in the next layer. MLP is fully-connected which means that each node of a layer is connected with every node in the next layer. The more layers used in an ANN's architecture, the "deeper" this architecture is considered.

Every node has an activation function and weights for every input and the aim of MLP is to train these weights from the input data, so the prediction is optimized. Without getting into much details, the algorithm for supervised regression problems has the main following steps: $(i)$ randomize the initial data, $(ii)$ train the network forwardly and predict the output, $(iii)$ calculate the error between the predicted and original values (usually with mean square error), $(iv)$ backpropagate the error (i.e. calculate it's derivative with respect to each weight) and update the model. This procedure is repeated for a number of iterations (epochs) until the error is minimized. 

One of the most popular activation function is ReLU (rectified linear unit) which is defined as
\begin{equation*}
    f(\mathrm{x}) = max(0,\mathrm{x})
\end{equation*}
where x is the input to a node and is a piecewise linear function. It outputs the input x directly if $\mathrm{x}>0$, otherwise, it will output zero. Another important parameter of the MLP algorithm, is the optimizer function that is used to minimize the error. Gradient decent-type optimizers are very common in neural networks, with Stochastic Gradient Decent (SGD) and Adaptive Moment Estimation (Adam) being very popular. 

Especially in big data problems, i.e. problems with large datasets, the complexity and nonlinear properties of the data may require a deeper architecture as described above, with more layers and nodes to successfully learn the features of that particular application. In general, using a much larger amount of nodes than the size of training samples, can lead to overfitting \cite{ZBMRV}. There are several methods developed to address this issue, with the two used in our particular applications being weight regularization and early stopping.

Weight regularization is a technique to reduce overfitting by adding a weight penalty term to the larger values of the weights of the network. There are two types: for L1 regularization the weight penalty is the sum of absolute values of the weights and for L2 regularization is the sum of squared values of the weights in the network. In both methods, the penalty term is multiplied by a small number $\lambda$, called regularization rate. Early stopping is a technique which aims to stop the training iteration, before the model begins to overfit. This is done by estimating the iteration after which the model's performance is not improved while the generalization error (starting learning from the noise of the data) is increasing. We refer to \cite{Yi} and the references therein, for a more detailed overview on how to avoid overfitting. 

\subsection{Hyperparameter tuning} 

As described above, all methods have different parameters which are essential for training and play a crucial role for the overall performance of each algorithm. These parameters are often referred as hyperparameters. Hyperparameter tuning is the process of optimizing the values or methods that define these parameters and can be a very challenging and time consuming task, depending on the complexity of each method.   

The most basic tuning method is grid search. With this method, a model is built for all possible combinations of the parameters, within a specific range of values. All models are then evaluated and the one with the best performance is selected. Tuning is usually performed together with k-fold cross-validation. The basic idea is to divide the dataset into k-unique parts, called folds, and use each of the k-times one fold as validate and the rest as train data. Then, evaluate the performance of each algorithm for all these couples. The final result is the  average of the folds. Application of cross-validation is important to reduce overfitting and ensure stable results. 
 
\subsection{Data pre-processing}

Pre-processing is an essential step in machine learning and refers to data cleaning and transformation from their raw forms to a more usable format, before building the model. The quality of the data used for training affects the ability of the models to learn. The scaling of input data in particular is very important. Scaling aims to standardize the samples in a specific scale or range. This can be useful e.g. in cases when our data varies in different orders of magnitude.

\begin{figure}[ht] 
   \begin{minipage}[b]{0.5\linewidth}
    \centering
    \includegraphics[scale=0.45]{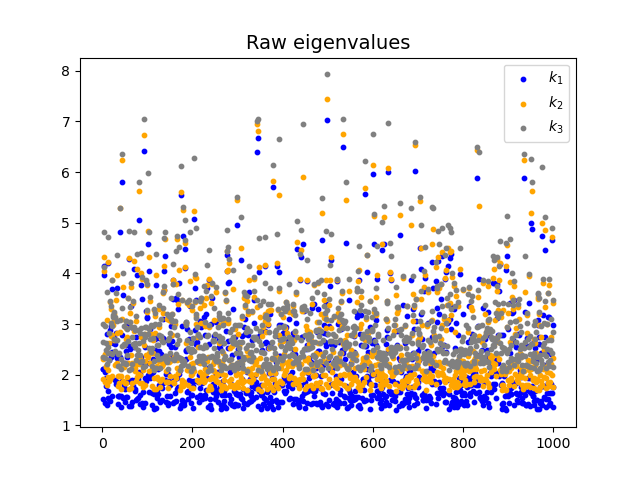} 
  \end{minipage}
  \begin{minipage}[b]{0.5\linewidth}
    \centering
    \includegraphics[scale=0.45]{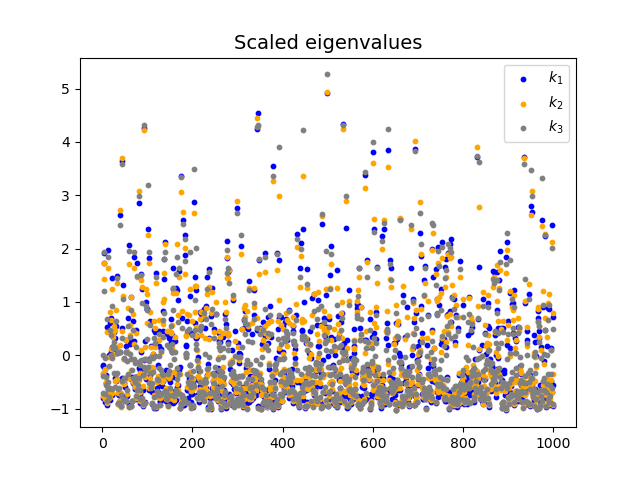} 
  \end{minipage} 
    \caption{Scaling eigenvalues with Standard Scaler.}
  \label{scale}
 \end{figure}

Here we use the Standard Scaler standardization, also known as Z-score normalization. With this method, data are scaled around their mean value, with unit standard deviation. Given a set of $\{\mathrm{X}_i\}_{i=1}^{N}$ values with mean value $\mu$ and standard deviation $\sigma$, the data are transformed to: 
\begin{equation*}
    \mathrm{z}_i=\frac{\mathrm{X}_i-\mu}{\sigma}.
\end{equation*}
As an example, figure \ref{scale} shows the distribution of 1000 samples for the first three transmission eigenvalues $k_1,\ k_2$ and $k_3$ of the  problem considered in \ref{sec421}. The eigenvalues before scaling have different order of magnitude (since the sequence of eigenvalues is increasing). After applying Standard Scaler, they are distributed around zero.  

In machine learning applications using Standard Scaler, both training and validating data must be standardized and the standardization should occur after the train-validate split to avoid data leakage. 

\subsection{Regression metrics}

Evaluation metrics are criteria for measuring the quality of a machine learning model. In the case of regression problems, some of the most common evaluation metrics are the mean square error (MSE), the root mean square error (RMSE) and the coefficient of determination $R^2$.  Assuming that we have a set of $\{\mathrm{y}_i\}_{i=1}^N$ predicted values and $\{\hat{\mathrm{y}_i}\}_{i=1}^N$ actual (ground truth) values, the MSE and the RMSE are given by:
\begin{equation*}
    \textrm{MSE} = \frac{1}{N}\sum_{i=1}^{N}(\mathrm{y}_i-\hat{\mathrm{y}_i})^2 \quad \text{and}\quad  \textrm{RMSE} = \sqrt{\frac{1}{N}\sum_{i=1}^{N}(\mathrm{y}_i-\hat{\mathrm{y}_i})^2}.
\end{equation*}
We note that RMSE has the same units as $\mathrm{y}$ and is in general more sensitive to the presence of outliers. 

Coefficient of determination $R^2$ is a good criterion on how well a model is trained, indicating the "closeness" of the predicted data to the actual. It is given by the following formula: 
\begin{equation*}
    R^2=1 - \frac{\mathrm{SS_{res}}}{\mathrm{SS_{tot}}} = 1- \frac{\sum_{i=1}^N (\mathrm{y}_i-\hat{\mathrm{y}_i})^2}{\sum_{i=1}^N (\mathrm{y}_i-\mu)^2},
\end{equation*}
where $\mathrm{SS_{res}}$ is the sum of squares of residuals, $\mathrm{SS_{tot}}$ the total sum of squares and $\mu$ the mean value of the $\mathrm{y}_i$. The $R^2$ varies from $(-\infty,1]$ and values near 1 indicate better predictions, closer to the actual values. Negative $R^2$ corresponds to cases when the prediction fits the model worse than a horizontal line would do.  Often, $R^2$ is presented in percentage $\%$ scale. 

\subsection{Feature evaluation}

When building a machine learning algorithm using input data X which consists of different features, we are often interested to use feature importance evaluation methods in order to rank the importance that each feature has on training the model. 
For example, in inverse problems like the ones we are studying here, we want to evaluate how important the eigenvalues are, according to their position in the eigenvalues sequence, in training the models we use. This could indicate which parameters play an important role on retrieving information about the unknown functions we want to determine. 

Here, we use an algorithm developed in \cite{HMFC} designed for deep neural networks, but is  appropriate for any supervised learning model. This method is called input perturbation algorithm and the main idea is to evaluate each trained model's efficiency with all of the features individually, after shuffling them. The most important features are those who produce a less accurate MSE score, when they are shuffled. 

\section{Numerical results for the inverse eigenvalue problems} \label{numer}

We firstly solve the direct spectral problems, to construct the datasets of eigenvalues which are used to train and validate the machine learning algorithms. Afterwards, we consider the inverse spectral problems where the unknown functions are reconstructed, as predictions of the underlying regression models. As described in section \ref{form}, for the inverse problems under consideration, there are uniqueness results supporting our study. Therefore, we have the theoretical background to proceed with solving these inverse problems using machine learning methods, since the knowledge of finite spectra is expected to provide information about the unknown functions \cite{LPR}.  

We note that all problems considered here, are multi-output i.e. supervised learning problems with multiple outputs to predict. Since all the output values (the coefficients of the differential operator) of size $m$ are related to the same input (the eigenvalues), we build a single model for predicting all $m$ outputs.

RF and kNN are implemented with the \textrm{scikit-learn} open source library \cite{skl} and MLP with the Keras neural network library, through the \textrm{TensorFlow} open source platform \cite{tens,ker}. We refer to \cite{He} for a guide on using Keras for deep neural networks. Furthermore, pre-processing, metrics and tuning of the algorithms is done using scikit-learn functions. Algorithms are written in \textrm{python} language. Computations are performed on a commercial intel core i7 computer, with 8GB RAM, and the run time for all models are in the order of a few minutes for training and a few hours for tuning. 

\subsection{Reconstructions of a symmetric potential} \label{invSL} 

To this end, we examine a simple example of the Sturm-Liouville problem (\ref{sl})-(\ref{slbc}), with Neumann boundary conditions and a symmetric potential, of the following form: 
\begin{eqnarray*}
     -y^{\prime\prime}+(1-e^{b(x-1/2)^2})y=\lambda y,&\quad 0<x<1,\\ 
     y^{\prime}(0)=y^{\prime}(1)=0,& \label{slsym}
\end{eqnarray*}
where $b\in\mathbb{R}$. We use \textsc{matslice}, to compute the $5$ first eigenvalues for a family of potentials with $-20\leq b \leq 20$ and step $0.04$. The potentials $q(x)=1-e^{b(x-1/2)^2}$ are discretized at $21$ points of $0\leq x \leq 1$, with step $0.05$. This scheme produces a set of $1001$ samples, where the training data consist of $5$ features \textrm{X} (the eigenvalues) and $21$ target values \textrm{y} (the discretized potentials). As a result, the multi-output models predict $m=21$ values. 

We randomly shuffle the data and split the $70\%$ for training and the $30\%$ for validating, which are identical for all models. The input samples are standardized using \textrm{Standard Scaler}. Afterwards, the inverse problem is considered, for a set of $5$ given eigenvalues corresponding to $4$ different potentials, for $b=-10\sqrt{2},\ -\sqrt{2}, \sqrt{2}$ and $10\sqrt{2}$. We note that these data are considered out-of-sample, in the sense that they are not included in the training/validating sets. 

Each machine learning model is tuned in their critical hyperparameters, using a grid search option and by performing cross-validation. For this specific training set the optimal parameters are: $(i)$ for RF the number of trees, is  "n$\_$estimators $= 100$", $(ii)$ for kNN, the number of neighbors is "n$\_$neighbors$=3$" and $(iii)$ for MLP we use an architecture of "$5-10-30-21$" with two hidden layers (i.e. 1041 trainable parameters), ReLU activation and Adam optimizer, L2 regularization with a regularization rate  $\lambda = 10^{-4}$ and early stopping. 

\begin{figure}[ht] 
   \begin{minipage}[b]{0.5\linewidth}
    \centering
    \includegraphics[scale=0.4]{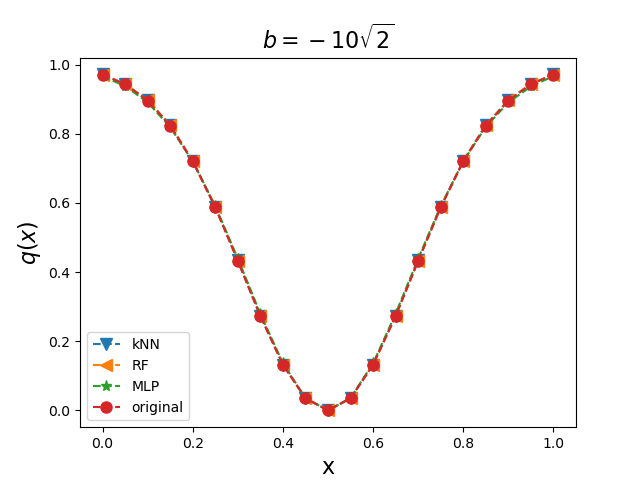} 
  \end{minipage}
  \begin{minipage}[b]{0.5\linewidth}
    \centering
    \includegraphics[scale=0.4]{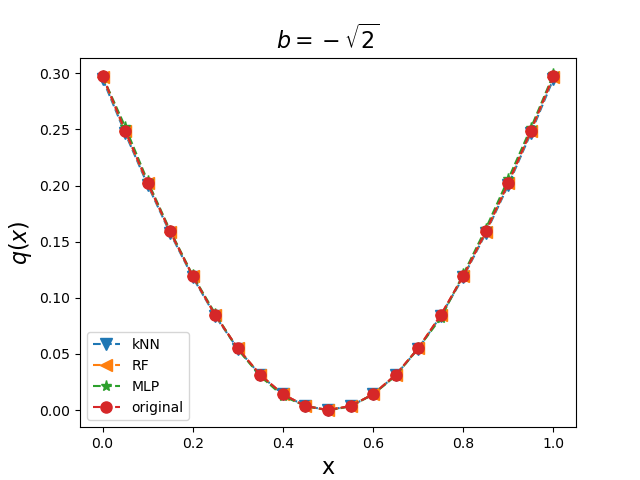} 
  \end{minipage} 
  \begin{minipage}[b]{0.5\linewidth}
    \centering
    \includegraphics[scale=0.4]{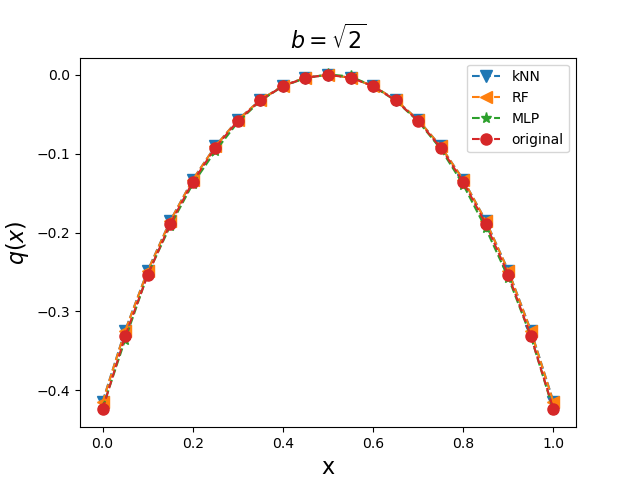} 
  \end{minipage}
  \begin{minipage}[b]{0.5\linewidth}
    \centering
    \includegraphics[scale=0.4]{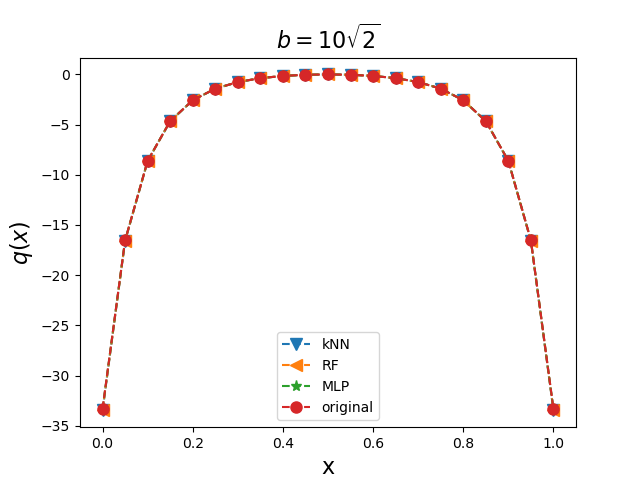} 
  \end{minipage} 
  \caption{Predictions of the symmetric potential $q(x)=1-
 e^{b(x-1/2)^2}$, using $5$ eigenvalues.}
  \label{slml}
\end{figure}

\begin{figure}[ht] 
  \centering
   \begin{minipage}[b]{0.5\linewidth}
    \centering
    \includegraphics[scale=0.4]{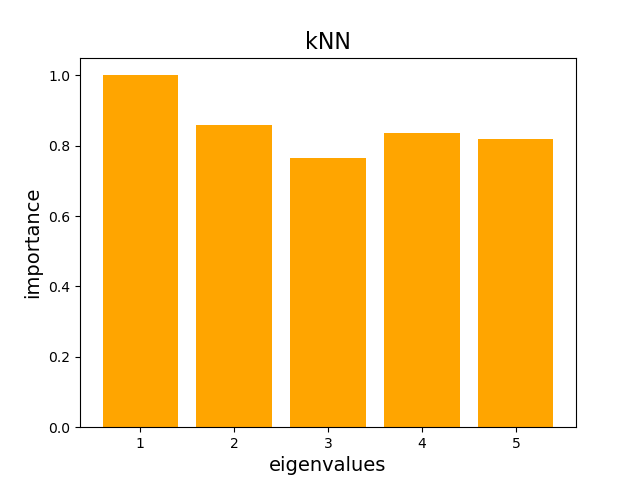} 
  \end{minipage}
  \begin{minipage}[b]{0.5\linewidth}
    \centering
    \includegraphics[scale=0.4]{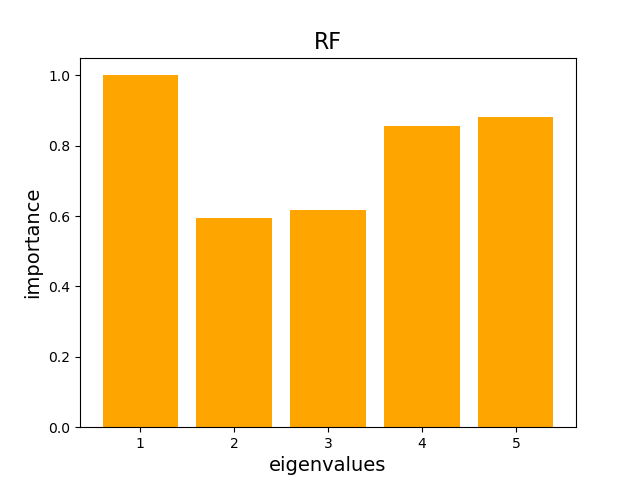} 
  \end{minipage} 
  \begin{minipage}[b]{0.5\linewidth}
    \centering
    \includegraphics[scale=0.4]{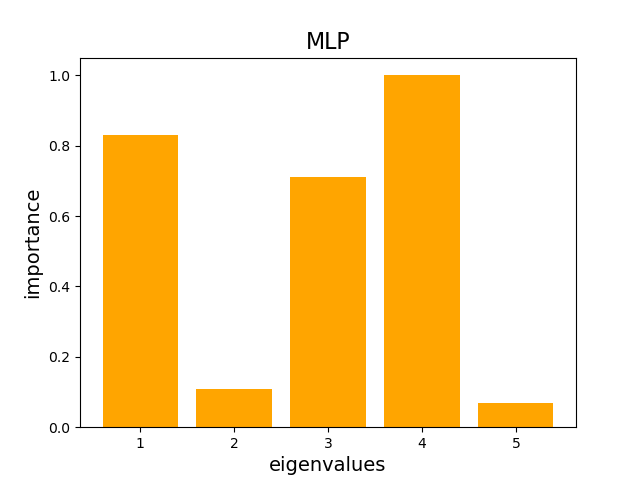} 
  \end{minipage}
  \caption{Eigenvalues importance according their position in the spectrum, for the Sturm-Liouville problem.}
  \label{rank_sl}
\end{figure}

We use $R^2$ regression score and RMSE , to compare the predictions with the corresponding ground truth potentials. For these $4$ testing potentials, the total $R^2$ scores are given in table \ref{r2sl}, showing great performance across all models. In table \ref{rmsesl} we present the RMSE for each potential, for all models under consideration. Figure \ref{slml} depicts for each model the predicted versus the corresponding original potentials. Finally, we use the input perturbation algorithm, to evaluate the importance of the features (i.e. the eigenvalues) for each model. As shown in figure \ref{rank_sl}, the lowest eigenvalue is important for training the models we used.  

\begin{table}
\centering
\caption{$R^2$ regression scores for training, validating and testing, for the inverse Sturm-Liouville eigenvalue problem.}
\addtolength{\tabcolsep}{5pt} 
\begin{tabular}{ |c||c||c||c|  }
 \hline
 model& train & validate & test\\
 \hline
   \hline
 kNN  & $99.99$    &$99.99$&   $99.99$\\
   \hline
  RF   & $99.99$    &$99.99$&   $99.99$\\
    \hline
   MLP  & $95.20$    &$95.20$&   $95.20$\\
 \hline
\end{tabular}\\
\label{r2sl}
  \end{table}

\begin{table}
\centering
\caption{RMSE between the predicted and the original potentials, for the inverse Sturm-Liouville eigenvalue problem.}
\addtolength{\tabcolsep}{5pt} 
\begin{tabular}{ |c||c||c||c|  }
 \hline
 potential& kNN & RF &MLP\\
 \hline
   \hline
 $b=-10\sqrt{2}$   & $0.038$    &$0.026$&   $0.069$\\
   \hline
  $b=-\sqrt{2}$   & $0.037$    &$0.019$&   $0.041$\\
    \hline
   $b=\sqrt{2}$   & $0.066$    &$0.064$&   $0.048$\\
     \hline
    $b=10\sqrt{2}$   & $0.229$    &$0.208$&   $0.12$\\
 \hline
 total&$0.115$&$0.104$&$0.072$\\
 \hline
\end{tabular}\\
\label{rmsesl}
  \end{table}

\subsection{Reconstructions of a piecewise constant refractive index} \
We now consider the transmission eigenvalue problem for a piecewise constant index, as described in subsection \ref{teppw}. We examine two different examples, one for refractive indices with two layers and one with four. 

In both examples, we use the Galerkin method to solve the direct problems and approximate the eigenvalues. We construct a basis with 30 eigenfunctions $\{\phi_i\}_{i=1}^{30}$ and compute the corresponding $30\times 30$ matrices $A^{(N)}, B^{(N)} $ and $C^{(N)}$ of the quadratic eigenvalue problem (\ref{a8}). The approximated transmission eigenvalues are calculated using the matlab function \textsc{polyeig} and are used for training/validating the machine learning algorithms. 

Furthermore, we compute the original transmission eigenvalues, as roots of the determinant of (\ref{detN}).  These eigenvalues are used for testing the algorithms and solving the inverse problems. We note that these original eigenvalues are considered out-of-sample of the training data and are assumed to to be "close enough" to the approximated eigenvalues, depending on the size $N$ of the Galerkin scheme. We choose to test the algorithms with these eigenvalues and not with a testing subset of the approximation samples X, because eigenvalues calculated from separation of variables are the original spectra associated with the scattering problem. Measurements coming from far-field data are expected to be close to these original eigenvalues.

For simplicity, we restrict ourselves in cases where $n(r)>2$ and use only real eigenvalues. It is shown that for refractive indices close to $1$, complex eigenvalues appear \cite{CCG}. By choosing larger refractive indices, complex eigenvalues with small modulus are avoided. We refer to \cite{GP,Pa} for the solution of the inverse transmission eigenvalue problem for piecewise constant indices, based on optimization methods. 

\subsubsection{Domain with two layers}\label{sec421}
Let $D$ be a disc of radius $r=1$ and the piecewise refractive index is $n=n_1$ for $0<r<d_1$ and $n=n_2$ for $d_1<r<1$. We solve the direct problem for $2\leq n_1,n_2 \leq 10$ with step $0.1$, for discontinuities at $0.1\leq d_1 \leq 0.9$ with step $0.1$. Thus, we construct a dataset with eigenvalues and refractive indices for all possible combinations of $n_1,n_2$ and $d_1$, which produce $59049$ samples. These training data consist of $6$ features \textrm{X} (the first $6$ eigenvalues) and $3$ target values \textrm{y} (the refractive index $(n_1,n_2,d_1)$). So the multi-output algorithms predict $m=3$ values simultaneously.  

In a similar fashion with the Sturm-Liouville problem, we shuffle the data, and for all models under consideration, we use the $70\%$ for training and the $30\%$ for validating. \textrm{Standard Scaler} is also used for scaling the input data set. For the inverse problem, we use the original first $6$ eigenvalues which are computed by separation of variables and test the machine learning algorithms, for $10$ different example refractive indices. 

After tuning, the optimal hyperparameters for each model are: $(i)$ for RF, the number of trees is "n$\_$estimators $= 200$", $(ii)$ for kNN, the number of neighbors is  "n$\_$neighbors$=3$" and $(iii)$ for MLP, we use an architecture of "$6-800-200-100-50-3$" with four hidden layers (i.e. $191103$ parameters) ReLU activation and Adam optimizer, L2 regularization with a $\lambda$ value of $10^{-4}$ and early stopping. 

\begin{figure}
\centering
  \begin{minipage}[b]{0.5\linewidth}
    \centering
    \includegraphics[scale=0.4]{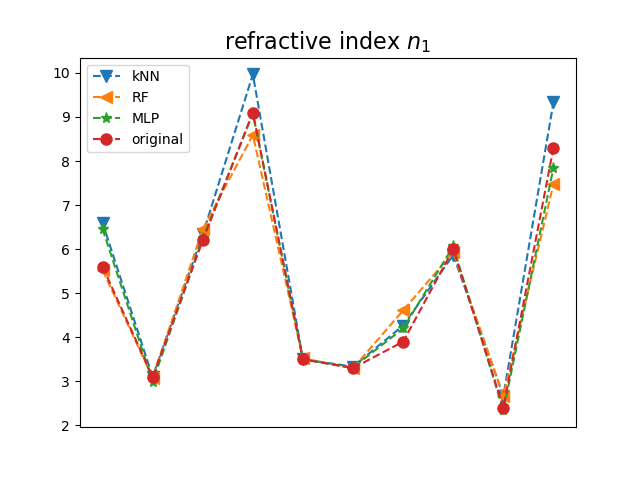} 
  \end{minipage}
  \begin{minipage}[b]{0.5\linewidth}
    \centering
    \includegraphics[scale=0.4]{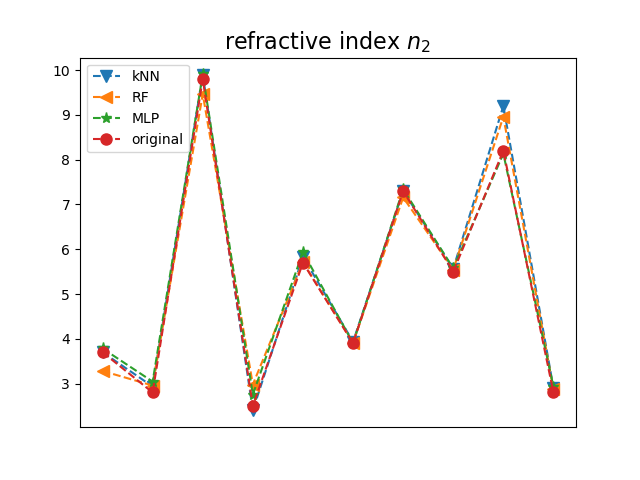} 
  \end{minipage} 
  \begin{minipage}[b]{0.5\linewidth}
    \centering
    \includegraphics[scale=0.4]{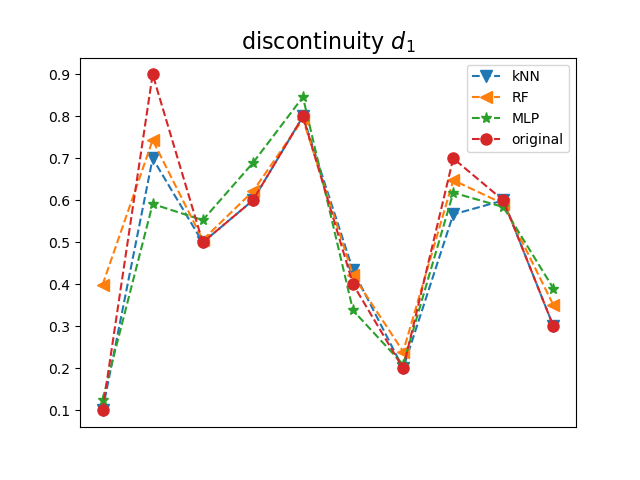} 
  \end{minipage}
  \caption{Predictions of the piecewise constant refractive index $(n_1,n_2,d_1)$.}
  \label{fig2Lml}
\end{figure}

The total $R^2$ scores for training, validating and testing can be found in table \ref{r22L}. 
\begin{table}
\centering
\caption{$R^2$ regression scores for training, validating and testing, for the inverse transmission eigenvalue problem of piecewise constant index with two layers.}
\addtolength{\tabcolsep}{5pt} 
\begin{tabular}{ |c||c||c||c|  }
 \hline
 model& train & validate & test\\
 \hline
   \hline
 kNN  & $96.27$    &$91.76$&   $93.99$\\
   \hline
  RF   & $98.88$    &$94.14$&   $91.54$\\
    \hline
  MLP  & $94.99$    &$94.74$&   $92.13$\\
 \hline
\end{tabular}\\
\label{r22L}
  \end{table}
Furthermore, table \ref{2Lml} includes the predictions for the inverse problem and in table \ref{rmse2L}  we compare the effectiveness of each model using the RMSE. The predictions for $n_1, n_2$ and $d_1$ respectively, are finally shown in figure \ref{fig2Lml}. 
\begin{table}
\centering
\caption{Predictions of piecewise constant index with two layers, from $6$ transmission eigenvalues.}
\addtolength{\tabcolsep}{5pt}    
\begin{tabular}{ |c||c||c||c|  }
\hline
original & \multirow{2}{*}{kNN}  & \multirow{2}{*}{RF}  & \multirow{2}{*}{MLP} \\
$(n_1,n_2,d_1)$ &  &   &  \\
  \hline
  \hline
 $5.6, 3.7, 0.1$   & $6.6,3.7,0.1$
& $5.53, 3.28, 0.4$
& $6.46,3.78,0.13$ \\
 \hline
  $3.1, 2.8, 0.9$   &$3.1,2.93,0.7$
& $3.08, 2.95, 0.74
$& $2.99,3.04,0.59$\\
 \hline
  $6.2, 9.8, 0.5$   & $6.33,9.9,0.5$ & $6.45, 9.47, 0.51
$ &$6.24,9.91,0.55$ \\
 \hline
  $9.1, 2.5, 0.6$   & $9.97,2.4,0.6
$&$8.59, 2.98, 0.62
$  &$9.08,2.75,0.69$\\
 \hline
  $3.5, 5.7, 0.8$   & $3.5,5.83,0.8
$& $3.52, 5.71, 0.79
$&$3.49,5.94,0.85$  \\
 \hline
  $3.3, 3.9, 0.4$   &$3.33,3.93,0.43
$ &$3.3,3.91,0.42
$  & $3.32,3.9,0.34$ \\
 \hline
  $3.9, 7.3, 0.2$   & $4.27,7.3,0.2
$ & $4.61, 7.16, 0.24
$& $4.2,7.34,0.21$ \\
 \hline
  $6, 5.5, 0.7$   & $5.87,5.57,0.57
$& $5.93, 5.54, 0.65
$& $6.06.5.57.0.62$
\\
 \hline
  $2.4, 8.2, 0.6$   &$2.67,9.2,0.6
$ & $2.67, 8.96, 0.59$
 & $2.34,8.14,0.58$ \\
 \hline
  $8.3, 2.8, 0.3$   & $9.33,2.9,0.3
$&$7.48, 2.87, 0.35
$& $7.85,2.93,0.39$  \\
\hline
\end{tabular}
\label{2Lml}
\end{table}  
\begin{table}
\centering
\caption{RMSE between the predicted and the original refractive indices, for the inverse transmission eigenvalue problem of piecewise constant index with two layers.}
\addtolength{\tabcolsep}{5pt}     
\begin{tabular}{ |c||c||c||c|  }
\hline
original & \multirow{2}{*}{kNN}  & \multirow{2}{*}{RF}  & \multirow{2}{*}{MLP} \\
$(n_1,n_2,d_1)$ &  &   &  \\
\hline
\hline
 $5.6, 3.7, 0.1$   & 0.76 & 0.548& 0.707   \\
 \hline
  $3.1, 2.8, 0.9$   & 0.373 & 0.357 &  0.485 \\
 \hline
  $6.2, 9.8, 0.5$   & 0.31 &  0.492 & 0.269 \\
 \hline
  $9.1, 2.5, 0.6$  &  0.71 & 0.635 & 0.395  \\
 \hline
  $3.5, 5.7, 0.8$   & 0.277 & 0.117 &  0.374 \\
 \hline
  $3.3, 3.9, 0.4$  &0.183& 0.117 &  0.192 \\
 \hline
  $3.9, 7.3, 0.2$  &0.46 & 0.647 &  0.418 \\
 \hline
  $6, 5.5, 0.7$   &  0.34 & 0.239 &  0.270 \\
 \hline
  $2.4, 8.2, 0.6$   & 0.773 & 0.683 & 0.226  \\
 \hline
  $8.3, 2.8, 0.3$   & 0.774 & 0.691 &   0.526 \\
 \hline
 total& 0.565 & 0.532 & 0.443 \\
 \hline
\end{tabular}
\label{rmse2L}
\end{table}

\begin{figure}[ht] 
  \centering
   \begin{minipage}[b]{0.5\linewidth}
    \centering
    \includegraphics[scale=0.4]{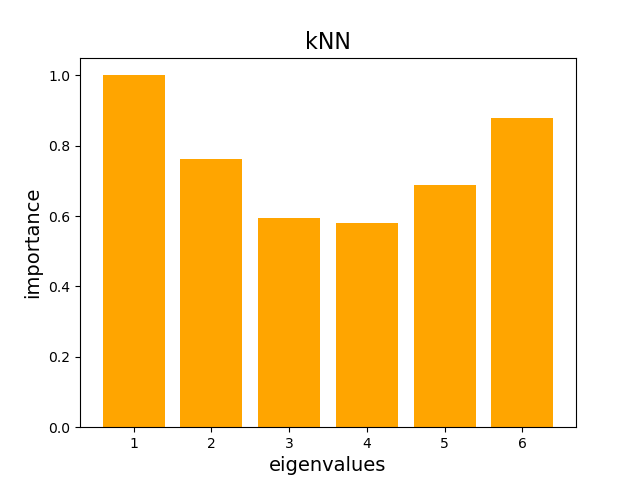} 
  \end{minipage}
  \begin{minipage}[b]{0.5\linewidth}
    \centering
    \includegraphics[scale=0.4]{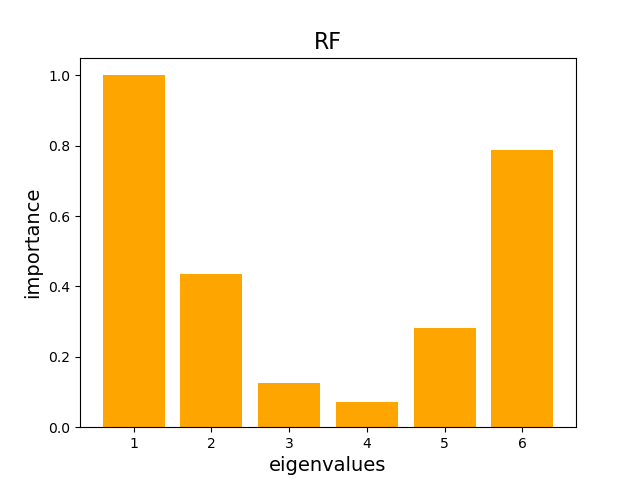} 
  \end{minipage} 
  \begin{minipage}[b]{0.5\linewidth}
    \centering
    \includegraphics[scale=0.4]{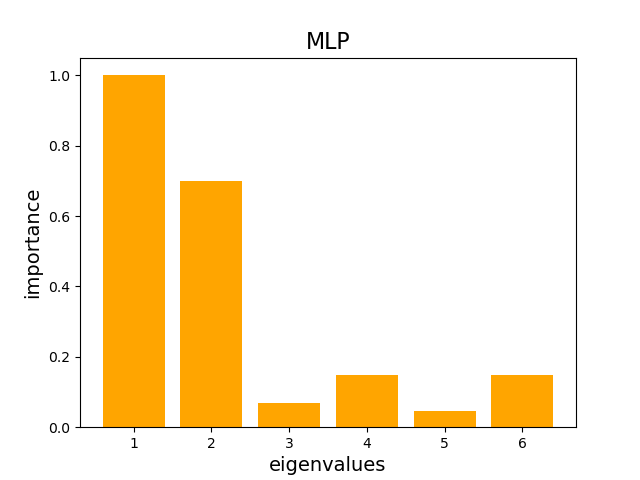} 
  \end{minipage}
  \caption{Eigenvalues importance according their position in the spectrum, for the transmission eigenvalue problem of piecewise constant index with two layers.}
  \label{rank_2L}
\end{figure}

We also estimate the importance of each eigenvalue as a training feature, using the input perturbation algorithm. In figure \ref{rank_2L} we notice that the first transmission eigenvalue is an important feature for all models used. The importance of the lowest real transmission eigenvalue in the inverse problem is also theoretically verified, since it is known that it can uniquely determine a constant refractive index of a medium with general shape \cite[Th. 3.1]{CCG}. 

Reconstructions are in good agreement with the original indices, showing an ability to train the models with the approximated eigenvalues and predict using the original transmission eigenvalues. We note that with these models we can recover the position of the discontinuity, which can be useful in applications like non-destructive testing of materials. 

\subsubsection{Domain with four layers} \label{inv4L} 
Finally, we examine the more complex case of a piecewise constant refractive index with four layers. The idea is similar to the two layer example discussed previously. We assume that each layer has width equal with $0.25$, i.e. $d_1=0.25,\ d_2=0.5$ and $d_3=0.75$. We solve the direct problem for all possible combinations of $2\leq n_1,n_2,n_3,n_4\leq 10$ with step $0.25$. This produces a dataset of $1185921$ samples, with $6$ eigenvalues as features \textrm{X} and $4$ target values \textrm{y}. Each multi-output algorithm predicts the $m=4$ outputs. 

Scaling, splitting and validating is done in the same way with the two layer case. Again, we use the first $6$ original eigenvalues of $10$ example refractive indices to test the predictive capability of each algorithm. 

The hyperparameters selected are: $(i)$ for RF, the number of trees is "n$\_$estimators $= 250$", $(ii)$ for kNN, the number of neighbors is  "n$\_$neighbors$=3$" and $(iii)$ for MLP, an architecture of "$6-800-500-500-200-50-4$" with five hidden layers (i.e. $767054$ parameters) is selected with SGD optimizer and ReLU activation, L2 regularization with a $\lambda$ value of $10^{-4}$ and early stopping. 

\begin{figure}
\centering
  \begin{minipage}[b]{0.5\linewidth}
    \centering
    \includegraphics[scale=0.4]{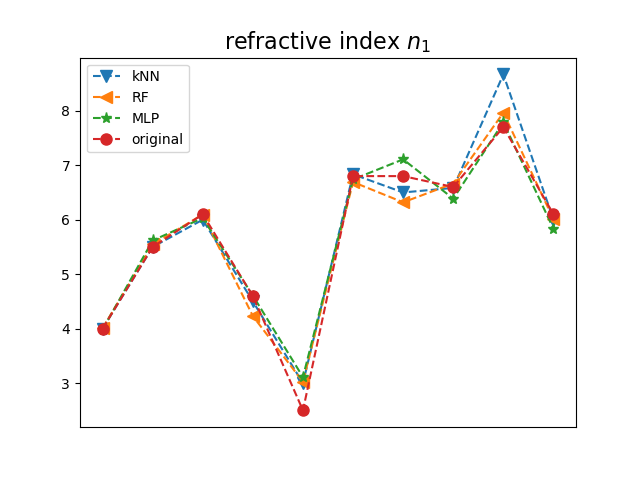} 
  \end{minipage}
  \begin{minipage}[b]{0.5\linewidth}
    \centering
    \includegraphics[scale=0.4]{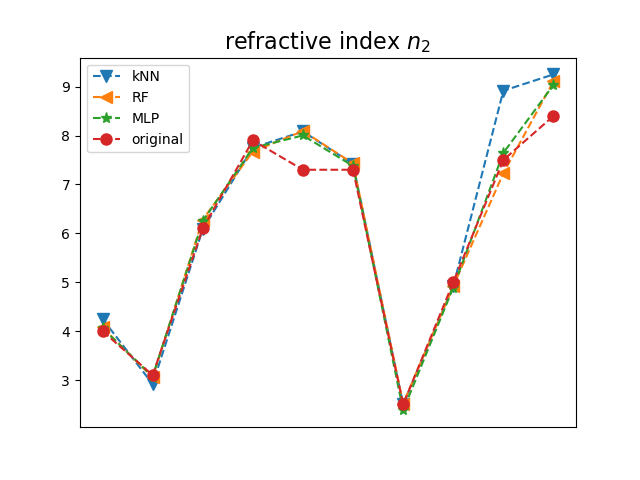} 
  \end{minipage} 
  \begin{minipage}[b]{0.5\linewidth}
    \centering
    \includegraphics[scale=0.4]{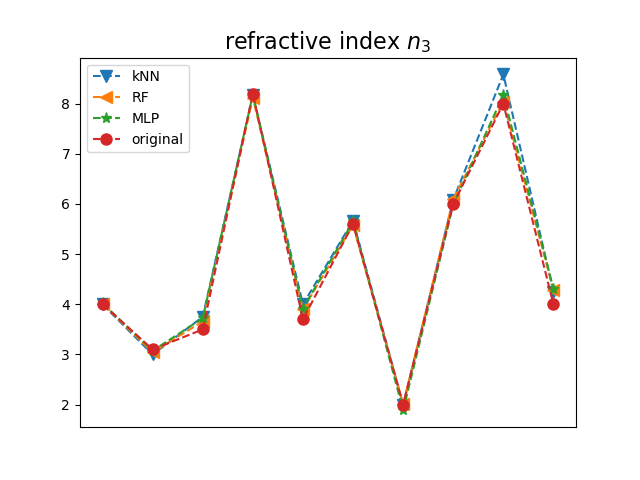} 
  \end{minipage}
   \begin{minipage}[b]{0.5\linewidth}
    \centering
    \includegraphics[scale=0.4]{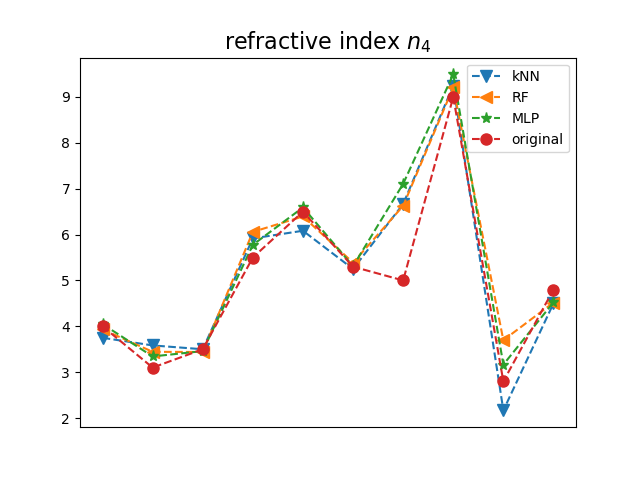} 
  \end{minipage}
  \caption{Predictions of the piecewise constant refractive index $(n_1,n_2,n_3,n_4)$.}
  \label{fig4Lml}
\end{figure}

\begin{figure}[ht] 
  \centering
   \begin{minipage}[b]{0.5\linewidth}
    \centering
    \includegraphics[scale=0.4]{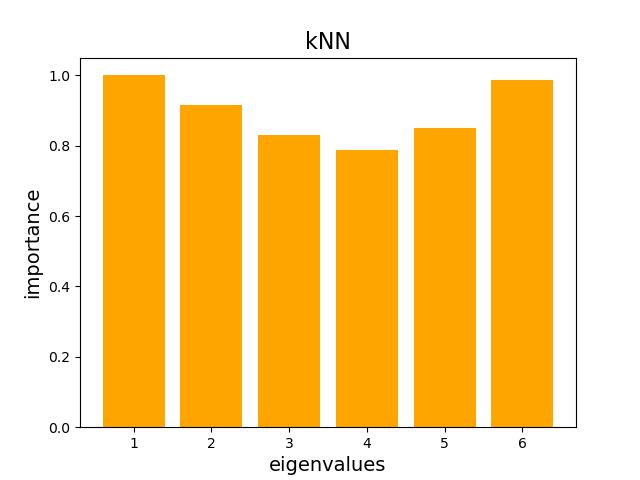} 
  \end{minipage}
  \begin{minipage}[b]{0.5\linewidth}
    \centering
    \includegraphics[scale=0.4]{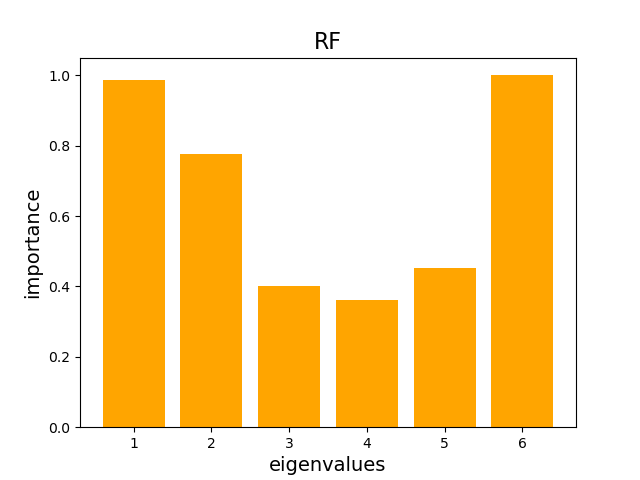} 
  \end{minipage} 
  \begin{minipage}[b]{0.5\linewidth}
    \centering
    \includegraphics[scale=0.4]{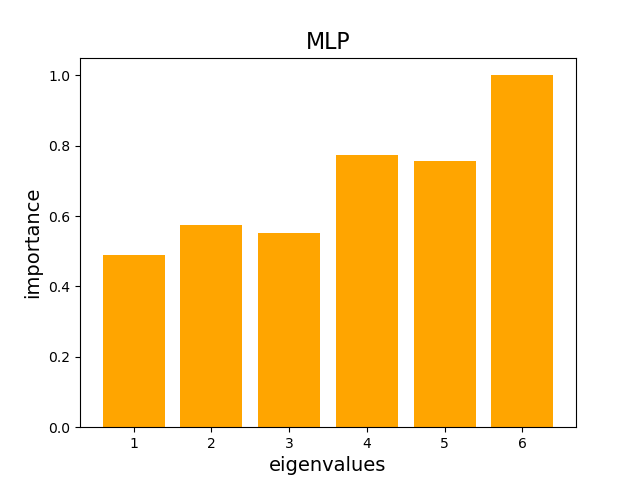} 
  \end{minipage}
  \caption{Eigenvalues importance according their position in the spectrum, for the transmission eigenvalue problem of piecewise constant index with four layers.}
  \label{rank_4L}
\end{figure}

We present the total $R^2$ scores in table \ref{r24L}. The reconstructions are included in table \ref{4Lml} and the corresponding RMSE in table \ref{rmse4L}. In figure \ref{fig4Lml} we present the per model predictions of each refractive index, while the evaluation of each eigenvalue is shown in figure \ref{rank_4L}. 
\begin{table}
\centering
\caption{$R^2$ regression scores for training, validating and testing, for the inverse transmission eigenvalue problem of piecewise constant index with four layers.}
\addtolength{\tabcolsep}{5pt} 
\begin{tabular}{ |c||c||c||c|  }
 \hline
 model& train & validate & test\\
 \hline
   \hline
 kNN  & $98.54$    &$96.99$&   $92.71$\\
   \hline
  RF   & $99.78$    &$98.43$&   $94.91$\\
    \hline
   MLP  & $99.79$    &$99.76$&   $94.38$\\
 \hline
\end{tabular}\\
\label{r24L}
  \end{table}

\begin{table}
\centering
\caption{Predictions of piecewise constant index with four layers, from $6$ transmission eigenvalues.}
\addtolength{\tabcolsep}{5pt}    
\begin{tabular}{ |c||c||c||c|  }
\hline
original & \multirow{2}{*}{kNN}  & \multirow{2}{*}{RF}  & \multirow{2}{*}{MLP} \\
$(n_1,n_2,n_3,n_4)$ &  &   &  \\
  \hline
  \hline
 $4,4,4,4$   & $4,4.25,4,3.75$ &$ 4.01,4.08,4.01,3.93$ & $3.98,4.06,4.02,4.05$ \\
 \hline
  $5.5,3.1,3.1,3.1$   &$5.5,2.92,3,3.58$&$ 5.55,3.05,3.04,3.44$ & $5.62,3.07,3.07,3.34$\\
 \hline
  $6.1,6.1,3.5,3.5$   & $6,6.08,3.75,3.5$ &$6.09,6.28,3.66,3.44 $ & $6.03,6.25,3.72,3.46$\\
 \hline
  $4.6,7.9,8.2,5.5$   & $4.5, 7.75, 8.17, 5.92$& $4.23,7.67,8.12,6.06 $& $4.61,7.75,8.18,5.78$\\
 \hline
  $2.5,7.3,3.7,6.5$   & $3,8.08,4,6.08$&$3.02,8.09,3.9,6.39$ & $3.11,8,3.91,6.59$ \\
 \hline
  $6.8,7.3,5.6,5.3$   &$6.83,7.42,5.56,5.25$&$6.69,7.43,5.58,5.39$ & $6.74,7.38,5.65,5.31$\\
 \hline
  $6.8,2.5,2,5$   & $6.5,2.5,2,6.67$ & $6.32,2.52,2,6.63$ & $7.11,2.38,1.89,7.11$\\
 \hline
  $6.6,5,6,9$   & $6.58,4.92,6.08,9.25$&$ 6.66,4.92,6.1,9.21$ & $6.39,4.89,5.99,9.49$\\
 \hline
  $7.7,7.5,8,2.8$   &$8.67,8.92,8.58,2.17$ & $ 7.96,7.23,8.04,3.7$& $7.79,7.64,8.17,3.17$ \\
 \hline
  $6.1,8.4,4,4.8$   & $6,9.25,4.25,4.5$& $6.01,9.11,4.29,4.52$& $5.83,9.03,4.31,4.52$\\
\hline
\end{tabular}
\label{4Lml}
\end{table}

\begin{table}
\centering
\caption{RMSE between the predicted and the original refractive indices, for the inverse transmission eigenvalue problem of piecewise constant index with four layers.}
\addtolength{\tabcolsep}{5pt}     
\begin{tabular}{ |c||c||c||c|  }
\hline
original & \multirow{2}{*}{kNN}  & \multirow{2}{*}{RF}  & \multirow{2}{*}{MLP} \\
$(n_1,n_2,n_3,n_4)$ &  &   &  \\
\hline
\hline
 $4,4,4,4$   & $0.42$& 0.225 & 0.201 \\
 \hline
  $5.5,3.1,3.1,3.1$   &$0.513$& 0.421 & 0.37 \\
 \hline
  $6.1,6.1,3.5,3.5$   & $0.367$ & 0.353 & 0.375\\
 \hline
  $4.6,7.9,8.2,5.5$   & $0.477$& 0.597 & 0.398\\
 \hline
  $2.5,7.3,3.7,6.5$   & $0.729$& 0.698 & 0.692 \\
 \hline
  $6.8,7.3,5.6,5.3$   & 0.271 & 0.309 & 0.237\\
 \hline
  $6.8,2.5,2,5$   & $0.92$ & 0.921 & 1.033\\
 \hline
  $6.6,5,6,9$   & $0.372$& 0.357 & 0.523\\
 \hline
  $7.7,7.5,8,2.8$   &$0.98$ & 0.698 & 0.467\\
 \hline
  $6.1,8.4,4,4.8$   & $0.686$& 0.641 & 0.634\\
 \hline
 total&$0.675$& 0.592 & 0.597 \\
 \hline
\end{tabular}
\label{rmse4L}
\end{table}

From the reconstructions of the first example in table \ref{4Lml} with constant refractive index $n=4$, we notice that all models can verify the homogeneity of the material, since the predictions are close to the nominal value. Furthermore, for the refractive index with only one discontinuity, i.e. $(5.5,3.1,3.1,3.1)$ and $(6.1,6.1,3.5,3.5)$, the predictions show that it's position is also correctly recovered. There are, however, cases where the predictions are not accurate enough, e.g. the index $(2.4, 7.4, 8.1, 2.1)$ is not successfully reconstructed. This is due to the fact that the six approximated eigenvalues, computed by the Galerkin scheme for this index are: $\{$2.8020,\ 3.0425,\ 3.1491,\ 3.4095,\ 3.7095,\ 4.0945$\}$ while the corresponding original are: $\{$2.5799,\ 2.6012,\ 2.9519,\ 3.3128, 3.6597,\ 4.0172$\}$, shown in figure \ref{4L_det}. This means that any inaccuracies that are within the training data, will eventually be passed to the predictions as well. This is a problem caused due to the size of approximation scheme, which is not large enough to achieve high accuracy in complex cases with big jumps between the layers of the index. We expect that a larger basis size $N$ could perform better approximations and thus more efficient training for the machine learning models. Also, we expect that using more eigenvalues than $6$ would perform better predictions, although this is also associated with the size $N$, since we have noticed that for larger eigenvalues the approximation error increases. 

Finally, from figure \ref{rank_4L}, we see that first and the sixth transmission eigenvalue are essential for training kNN and RF. Sixth eigenvalue is the most important feature for MLP.

\subsection{Discussion} 

In all examples examined, we achieved satisfactory predictions, which show that the models considered in this paper, are capable to numerically solve these inverse spectral problems. However,  
predictions with out-of-sample data can be considered reliable only for eigenvalues that fall within the range of their training data.

For the simple case of the inverse Sturm-Liouville eigenvalue problem in \ref{invSL}, we expected that the machine learning models would perform pretty well, due to the self-adjointness and the continuous dependence between the spectral data and the potentials \cite{LPR, RS}. 

The inverse transmission eigenvalue problem, had a more challenging solution due to it's non-self-adjoint and more complex nature. Models for the problem considered in \ref{inv4L} with four layers, can (surprisingly) achieve more  successful results, compared to those of the two layer problem in \ref{sec421}. This is mainly due to the fact that the discontinuities are fixed, while the discontinuity in \ref{sec421} is variable. Our results show that the position of the discontinuity is in general more difficult to be predicted and this is responsible for the relatively lower $R^2$ scores in table \ref{r22L} compared to table \ref{r24L}. To study more complex problems, consisting of multiple layers and variable discontinuity positions, it would require the production of very large scale training data as well as a lot more computational time and power to train the algorithms. It would also need more input eigenvalues and a larger size for the approximation scheme. Nevertheless, the overall results show that the algorithms are able to get trained by approximated data and predict using original, out-of-sample data. This can be useful in applications like non-destructive testing. The above ideas could, furthermore, be applied in general domains without spherical symmetries. By utilizing the spectral-Galerkin method the, required for training, spectral data can be calculated using (\ref{a8}). Although, for geometries that do not allow the separation of variables technique, one should test the machine learning algorithms with original transmission eigenvalues which are measured from far-field data \cite{CCH2,CCM}.

The feature evaluation algorithm, gave us indications about the importance of each eigenvalue in training the models. 
Further investigation of the effect of each eigenvalue to the inverse problem, both from theoretical and numerical point of view, could be the subject of a future study.

Finally, we would like to mention that for all examples studied in section \ref{numer}, we also performed a 10-time random split to the input sets and averaged the output predictions. The outcome showed stable $R^2$ scores, for all train-validate-test datasets. This was done to ensure the stability of the samples used as well as the performance of each model.

\section{Summary}
In this paper, we studied the application of supervised machine learning regression models to the numerical solution of inverse eigenvalue problems. We considered the self-adjoint case of the classic Sturm-Liouville eigenvalue problem with symmetric potentials as well as the non-self-adjoint transmission eigenvalue problem for a spherically symmetric piecewise constant index. We solved the direct spectral problems to produce datasets of eigenvalues for given potentials and refractive indices, which where then used to train and validate the machine learning algorithms. 
All models achieved successful predictions. The efficiency of each model was compared using $R^2$ regression score and RMSE. Finally, we evaluated the importance of the eigenvalues, according to their position in the spectrum, by implementing an input perturbation algorithm. The above machine learning models (and other potential regression methods) could be tested in other inverse eigenvalue problems, provided sufficient samples of input spectral data are obtained to successfully train the algorithms. 

\subsection*{Data availability}  The data that support this study are available by the authors, upon reasonable request.

\subsection*{Conflict of interest} All authors declare that they have no conflicts of interest. 

\subsection*{Acknowledgements} 
Authors would like to thank the n\&k Technology, Inc. engineering team, for the fruitful conversations about applications of machine learning.

\end{document}